
\documentstyle{amsppt}
\baselineskip18pt
\magnification=\magstep1
\pagewidth{30pc}
\pageheight{45pc}

\hyphenation{co-deter-min-ant co-deter-min-ants pa-ra-met-rised
pre-print pro-pa-gat-ing pro-pa-gate
fel-low-ship Cox-et-er dis-trib-ut-ive}
\def\leaderfill{\leaders\hbox to 1em{\hss.\hss}\hfill}
\def\A{{\Cal A}}

\def\H{{\Cal H}}

\def\L{{\Cal L}}

\def\Pl{{\Cal P}}
\def\Sy{{\Cal S}}\

\def\ldescent#1{{\Cal L (#1)}}
\def\rdescent#1{{\Cal R (#1)}}

\def\afn{{\text {\bf a}}}
\def\tr{{\text {\rm tr}}}

\def\d{{\delta}}

\def\th{{\theta}}

\def\l{{\lambda}}

\def\t{{\tau}}

\def\T{{\widetilde T}}
\def\te{\widetilde t}

\def\tmu{\tilde{\mu}}

\def\b0{\text{\bf 0}}

\def\ra{{\ \longrightarrow \ }}

\def\lan{{\langle}}
\def\ran{{\rangle}}

\def\real{{\Bbb R}}

\def\zed{{\Bbb Z}}
\def\kyu{{\Bbb Q}}
\def\enn{{\Bbb N}}

\def\Im{\text{\rm Im}}

\def\boxit#1{\vbox{\hrule\hbox{\vrule \kern3pt
\vbox{\kern3pt\hbox{#1}\kern3pt}\kern3pt\vrule}\hrule}}
\def\rabbit{\vbox{\hbox{\kern0pt
\vbox{\kern0pt{\hbox{---}}\kern3.5pt}}}}

\def\tableau#1{
        \hbox {
                \hskip -10pt plus0pt minus0pt
                \raise\baselineskip\hbox{
                \offinterlineskip
                \hbox{#1}}
                \hskip0.25em
        }
}

\def\tabCol#1{
\hbox{\vtop{\hrule
\halign{\strut\vrule\hskip0.5em##\hskip0.5em\hfill\vrule\cr\lower0pt
\hbox\bgroup$#1$\egroup \cr}
\hrule
} } \hskip -10.5pt plus0pt minus0pt}

\def\CR{
        $\egroup\cr
        \noalign{\hrule}
        \lower0pt\hbox\bgroup$
}



\def\blank#1#2{
\hbox to #1{\hfill \vbox to #2{\vfill}}
}


\def\strut{\vrule height10pt depth5pt width0pt}

\def\secz{1}
\def\secy{2}
\def\seca{3}
\def\secb{4}
\def\secc{5}
\def\secd{6}
\def\sece{7}
\def\secf{8}
\def\secg{9}
\def\qvk{K}
\def\tlqvk{TL_{\qvk}(E_n)}
\def\tlqvkd{TL_{\qvk'}(E_n)}
\def\almd{(1+v^{-2})}

\topmatter
\title On the Markov trace for Temperley--Lieb algebras of type $E_n$
\endtitle

\author R.M. Green \endauthor
\affil Department of Mathematics \\ University of Colorado \\
Campus Box 395 \\ Boulder, CO  80309-0395 \\ USA \\ {\it  E-mail:}
rmg\@euclid.colorado.edu \\
\newline
\endaffil

\abstract 
We show that there is a unique Markov trace on the tower
of Temperley--Lieb type quotients of Hecke algebras of Coxeter type $E_n$
(for all $n \geq 6$).  We explain in detail how this trace may be computed
easily using tom Dieck's calculus of diagrams.  As applications, we show
how to use the trace to show that the diagram representation is faithful,
and to compute leading coefficients of certain Kazhdan--Lusztig polynomials.
\endabstract

\subjclass 20C08, 20F55, 57M15 \endsubjclass

\endtopmatter


\head \secz. Introduction \endhead

In the paper \cite{{\bf 17}}, Jones introduced a certain Markov trace on 
the tower of Hecke algebras $\H(A_{n-1})$ associated to the Coxeter groups 
$\Sy_n = W(A_{n-1})$, which are the symmetric groups.
When Jones' trace is restricted to one of the algebras $\H = \H(A_{n-1})$, 
it is degenerate, but its radical is an ideal, $J$, of $\H$ and so we
obtain a generically nondegenerate trace on the algebra $\H/J$, which is the
Temperley--Lieb algebra $TL_n$ occurring in statistical mechanics \cite{{\bf 25}}
(the trace is the matrix trace of a transfer matrix algebra).

In \cite{{\bf 19}}, Kazhdan and Lusztig introduced a remarkable polynomial
$P_{x, w}(q)$ for any elements $x, w$ in a Coxeter group $W$.  These 
polynomials have important applications in representation theory.
Although the polynomials have an elementary definition, the only
obvious way to compute them is using a rather complicated recurrence
relation.  One of the main obstructions to computing the polynomials 
efficiently is a fast way to compute the integer $\mu(x, w)$, which is the
coefficient of $q^{(\ell(w) - \ell(x) - 1)/2}$ in $P_{x, w}(q)$.
In \cite{{\bf 12}}, the author showed how Jones' trace can be used to compute
the leading coefficients $\mu(x, w) \in \zed$ in the case where $x$ and $w$
are fully commutative elements of $W$ (in the sense of \cite{{\bf 24}}).  
In this paper, we will
investigate the analogous phenomenon in Coxeter type $E_n$.
This includes Coxeter groups of types $A$ and $D$ as special cases.

The algebras $TL_n$ may be defined in terms of generators and relations
in a way that generalizes readily to Coxeter systems of other types.
These generalized Temperley--Lieb algebras have been studied for Coxeter 
type $E_n$ by a number of people \cite{{\bf 2}, {\bf 3}, {\bf 7}}.
Although the Coxeter groups of type $E_n$ are infinite for $n > 8$, the 
Hecke algebra quotient $TL(E_n)$ in this case is still finite dimensional.
In \cite{{\bf 2}}, tom Dieck constructed a diagrammatic representation
of $TL(E_n)$, although the question of whether this is a realisation---a 
faithful representation---is not tackled.  In \S\secg, we will prove

\proclaim{Theorem \secz.1}
The diagrammatic representation of $TL(E_n)$ given in \cite{{\bf 2}} is injective.
\endproclaim

The closing remarks of \cite{{\bf 2}} state
without proof that this representation can be used to define a Markov trace
on the tower of algebras $TL(E_n)$.  In Theorem \secf.11, we will prove
this claim and furthermore we will show that there is a unique such Markov
trace.  Although this is similar to what happens in type $A$, 
the analogous claim for Coxeter type $D$ is false.

This trace is also remarkable for other reasons: after suitable rescaling,
it is a tabular trace in the
sense of \cite{{\bf 10}}, and a generalized Jones trace in the sense of \cite{{\bf 12}}.
The fact that the trace is tabular implies that it is (generically) 
nondegenerate on the algebras $TL(E_n)$.  The fact that we have a generalized
Jones trace will lead to the following theorem (proved in \S\secg)
where the monomial basis elements $b_w$ are defined in \S\seca.

\proclaim{Theorem \secz.2}
Let $\{b_w : w \in W_c\}$ be the monomial basis of $TL(E_n)$ indexed by the
fully commutative Coxeter group elements, and let $\tr$ be the unique Markov
trace on the tower of algebras $TL(E_n)$.  If $x, y \in W_c$, then the
coefficient of $v^{-1}$ in $\tr(b_x b_{y^{-1}})$ (after expansion as a power
series) is $\tmu(x, y)$, where $$
\tmu(x, y) = \cases
\mu(x, y) & \text{ if } x \leq y,\cr
\mu(y, x) & \text{ if } x \not\leq  y,\cr
\endcases
$$ and $\mu(a, b)$ is the integer defined in \cite{{\bf 19}}.
\endproclaim

We will also show in \S\secg\  how $\tmu(x, y)$ may be evaluated 
non-recursively using the diagram calculus.

\head \secy. Traces and Markov traces \endhead

By a {\it trace} on an $R$-algebra $A$, we mean an $R$-linear map $t : A \ra R$
such that $t(ab) = t(ba)$ for all $a, b \in A$.  The {\it radical} of the
trace is the set of all $a \in A$ such that $t(ab) = 0$ for all $b \in A$.
The radical is always an ideal of $A$, and if it is trivial, the trace is
said to be {\it nondegenerate}.  In any case, if $I$ is the radical of $t$,
then $t$ induces a nondegenerate trace on the quotient algebra $R/I$.

The set of traces on an $R$-algebra $A$ has a natural $R$-module structure.
In the special case where $\rho$ is a representation of an $R$-algebra $A$,
then the matrix trace associated to $\rho$ is a trace in the above sense,
which means that, if $A$ is semisimple, the Grothendieck group of $A$ gives 
a $\zed$-lattice in the space of traces, generated by the traces of the
simple modules.

We will be particularly concerned with algebras where the base ring $R$ is 
obtained by extending scalars from the ring of Laurent polynomials $\A = 
\zed[v, v^{-1}]$ to some ring $F \otimes \A$.  This has the effect of 
specializing the parameter $v$ to an invertible element of $F$.
In this situation, a trace is called {\it generically nondegenerate} if it is 
nondegenerate as a trace over $\A$, and if it also remains nondegenerate as
a trace over $F \otimes \A$ for all but finitely many specializations of
$v$.

Suppose now that $R$ is an integral domain and 
$\{A_n : n \geq N\}$ is a family of unital $R$-algebras such that
$A_n$ is a subalgebra of $A_{n+1}$ for all $n \geq N$.  Let $A_{\infty}$
be the associated direct limit.  Suppose also that
there is a set of elements $\{g_n : n \in \enn\}$ such that $g_{n+1} \in 
A_{n+1} \backslash A_n$ for all $n$ and such that
$\{g_n : n \leq M\}$ is an algebra generating set for $A_M$.  Following
\cite{{\bf 5}, \S4}, we may now introduce the notion of Markov trace.

\definition{Definition \secy.1}
Maintain the above notation, and let $F$ be a field containing $R$.
A {\it Markov trace} on $A_\infty$ with parameter $z \in F$
is an $F$-linear map 
$\t : A_\infty \ra F$ 
satisfying the following conditions:
\item{\rm (i)}{$\t(1) = 1$;}
\item{\rm (ii)}{$\t(h b_{n+1}) = z \t(h)$ for $n \geq N$ and $h \in A_n$;}
\item{\rm (iii)}{$\t(hh') = \t(h'h)$ for all $h, h' \in A_\infty$.}
\enddefinition

Jones \cite{{\bf 17}} proved that there is a unique Markov
trace with parameter $z$ on the tower of Hecke algebras of type $A_n$, and 
that the only one of these traces that passes
to the Temperley--Lieb quotient is the one with parameter 
$z = (v + v^{-1})^{-1}$.  This is an important observation in the construction
of the Jones polynomial, because conditions (ii) and (iii) for the trace are
what is needed to ensure that the polynomial is invariant under the two
types of ``Markov move''.

Some other notable work on Markov traces includes that of Geck and 
Lambropoulou \cite{{\bf 4}}, who classified the Markov traces in Coxeter 
types $B$ and $D$, using a suitable extension of the above definition.
Lambropoulou \cite{{\bf 20}} extended this work (in type $B$) to generalized and
cyclotomic Hecke algebras of type $B$.

For the purposes of studying Temperley--Lieb type quotients of Hecke
algebras, a better definition of Markov traces seems to be one that 
appears in work of Seifert \cite{{\bf 22}} and
recent work of Gomi \cite{{\bf 6}, Definition 3.7}.  In this case, one 
retains conditions (i) and (iii) 
of Definition \secy.1 and replaces condition (ii) by the requirement that $$
\tau(a T_s) = z_s T(a)
$$ whenever we have $a \in {\Cal H}(W_I)$ for some parabolic subgroup
$W_I$ corresponding to $I \subseteq S \backslash \{s\}$.  (In other words,
we require condition (ii) to hold for all generators of $A_{n+1}$, not just
one particular generator.)  Here, $z_s$ is
an indeterminate depending on the conjugacy class of $s$ in $W$.  

In this paper, we will restrict our attention to the tower of algebras
$TL(E_n)$, and in this case, the above definitions happen to agree; however,
they do not agree in the corresponding question for type $D_n$.  In the latter
case, it can be shown that the Seifert--Gomi formulation produces a unique 
Markov trace, and Definition \secy.1 does not.

\head \seca. The algebras $TL(E_n)$ \endhead

Let $X = X(E_n)$ be a Coxeter graph of type $E_n$, where $n \geq 6$.  
Following \cite{{\bf 3}}, we label the vertices of $X$ by
$0, 1, \ldots, n-1$ in such a way that $1, 2, 3, \ldots, n-1$ lie in a straight
line, and such that $3$ is the unique vertex of degree $3$, which is adjacent
to $2$, $4$ and $0$.  Figure 1 shows the case $n = 6$.

\topcaption{Figure 1} Coxeter graph of type $E_6$ \endcaption
\centerline{
\hbox to 2.138in{
\vbox to 0.819in{\vfill
        \includegraphics{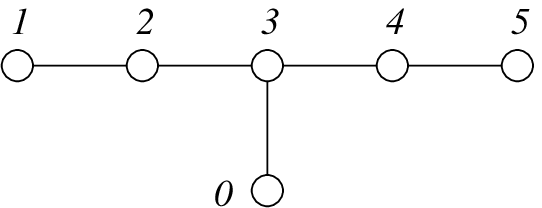}
}
\hfill}
}

Let $W(E_n)$ be the associated Coxeter group with distinguished
set of generating involutions $$
S(E_n) = \{s_i : i \text{ is a vertex of } X(E_n)\}
.$$  In other words, $W = W(E_n)$ is 
given by the presentation $$
W = \lan S(E_n) \ | \ (st)^{m(s, t)} = 1 \text{ for } m(s, t) < \infty \ran
,$$ where $m(s, s) = 1$, $m(s, t) = 2$ if $s$ and $t$ are not adjacent in $X$,
and $m(s, t) = 3$ if $s$ and $t$ are adjacent in $X$.
The elements of $S = S(E_n)$ are distinct as group elements, and $m(s, t)$ 
is the order of $st$.  Denote by $\H_q = \H_q(E_n)$ the Hecke
algebra associated to $W$.  This is a $\zed[q, q^{-1}]$-algebra
with a basis consisting of (invertible) elements $T_w$, with $w$ ranging over 
$W$, satisfying $$T_s T_w = 
\cases
T_{sw} & \text{ if } \ell(sw) > \ell(w),\cr
q T_{sw} + (q-1) T_w & \text{ if } \ell(sw) < \ell(w),\cr
\endcases$$ where $\ell$ is the length function on the Coxeter group
$W$, $w \in W$, and $s \in S$.  If $n > 8$, the group $W$ is infinite and
$\H_q$ has infinite rank as an $\A$-algebra.

For the applications we have in mind, it is convenient to 
extend the scalars of $\H_q$ to produce an $\A$-algebra $\H$,
where $\A = \zed[v, v^{-1}]$ and $v^2 = q$, and to define a scaled version
of the $T$-basis, $\{\T_w : w \in W\}$, where $\T_w := v^{-\ell(w)} T_w$.
We will write $\A^+$ and $\A^-$ for $\zed[v]$ and $\zed[v^{-1}]$, respectively.

A product $w_1w_2\cdots w_n$ of elements $w_i\in W$ is called
{\it reduced} if \newline $\ell(w_1w_2\cdots w_n)=\sum_i\ell(w_i)$.  We reserve
the terminology {\it reduced expression} for reduced products 
$w_1w_2\cdots w_n$ in which every $w_i \in S$.  We write $$
\ldescent{w} = \{s \in S : \ell(sw) < \ell(w)\}
$$ and $$
\rdescent{w} = \{s \in S : \ell(ws) < \ell(w)\}
.$$  The set $\ldescent{w}$ (respectively, $\rdescent{w}$) is called the 
{\it left} (respectively, {\it right}) {\it descent set} of $w$.

Call an element $w \in W$ {\it complex} if it can be written 
as a reduced product $x_1 w_{ss'} x_2$, where $x_1, x_2 \in W$ and
$w_{ss'}$ is the longest element of some rank 2 parabolic subgroup 
$\lan s, s'\ran$ such that $s$ and $s'$ correspond to adjacent vertices
in the Coxeter graph $E_n$.
Denote by $W_c(E_n)$ the set of all elements of $W$
that are not complex.  The elements of $W_c = W_c(E_n)$ are the 
{\it fully commutative} elements of \cite{{\bf 24}}; they are characterized by 
the property that any two of their reduced expressions may be obtained 
from each other by repeated commutation of adjacent generators.

Let $J(E_n)$ be the two-sided ideal of $\H$ generated by the elements $$
T_1 + T_s + T_t + T_{st} + T_{ts} + T_{sts}
,$$ where $(s, t)$ runs over all pairs of elements of $S$ for which
$m(s, t) = 3$.  Following Graham \cite{{\bf 7}, Definition 6.1}, we define the 
{\it generalized Temperley--Lieb algebra} $TL(E_n)$ to be
the quotient $\A$-algebra $\H(E_n)/J(E_n)$.  We denote the corresponding
epimorphism of algebras by $\th : \H(E_n) \ra TL(E_n)$.  
Let $t_w$ (respectively, $\te_w$) denote the image in 
$TL(E_n)$ of the basis element $T_w$ (respectively, $\T_w$) of $\H$.
If $s \in S$, we define $b_s \in TL(E_n)$ by $b_s = v^{-1} 1 + \te_s$.

A more convenient description of $TL(E_n)$ for the purposes of this paper is 
by generators and relations (as in \cite{{\bf 3}, \S2.2}).  Since the Laurent
polynomial $v + v^{-1}$ occurs frequently, we denote it by $\d$.

\proclaim{Proposition \seca.1}
As a unital $\A$-algebra, $TL(E_n)$ is given by generators $\{b_s : s \in S\}$
and relations $$\eqalign{
b_s^2 &= \d b_s, \cr
b_s b_t &= b_t b_s \text{\quad if } m(s, t) = 2, \cr
b_s b_t b_s &= b_s \text{\quad if } m(s, t) = 3. \cr
}$$ 
\qed\endproclaim

The following basis theorem will be used freely in the sequel.

\proclaim{Theorem \seca.2 \cite{{\bf 3}, {\bf 7}}}
\item{\rm (i)}
{The set $\{\te_w : w \in W_c\}$ is a free $\A$-basis for $TL(E_n)$.}
\item{\rm (ii)}{If 
$w \in W_c$ and $w = s_{i_1} s_{i_2} \cdots s_{i_r}$ is reduced, then the
element $$
b_w = b_{s_{i_1}} b_{s_{i_2}} \cdots b_{s_{i_r}}
$$ is a well-defined element of $TL(E_n)$.}
\item{\rm (iii)}
{The set $\{b_w : w \in W_c\}$ is a free $\A$-basis for $TL(E_n)$.}
\endproclaim

\demo{Proof}
Part (i) is due to Graham \cite{{\bf 7}, Theorem 6.2}.  Parts (ii) and (iii) are
stated by Fan in \cite{{\bf 3}, \S2.2}, and more details may be found in 
\cite{{\bf 13}, Proposition 2.4}.
\qed\enddemo

\definition{Definition \seca.3 \cite{{\bf 3}, \S2.3}}
Let $P = P(n)$ denote the set of subsets of the Coxeter graph $E_n$ that 
consist of 
non-adjacent vertices.  We allow $P$ to include the empty set, $\emptyset$. 
For any $A \in P$, let $i(A)$
be the product of the elements of $S(E_n)$ corresponding to the vertices in $A$
(with $i(\emptyset) = 1$);
note that the order of the product is immaterial since the vertices in $A$ 
correspond to commuting generators.
Let $A, B \in P$.  We say that $A$ and $B$ are neighbours if and only if 
$1 + \#(A \cap B) = \#A = \#B$, and the two vertices in $(A \cup B) \backslash
(A \cap B)$ are adjacent in $E_n$. Define an
equivalence relation on $P$ by taking the reflexive and transitive closure 
of the relation $A \sim B$ if $A$ and $B$ are neighbours. Let $\bar{P}$ 
denote the set $P/\sim$ .
\enddefinition

\example{Example \seca.4}
In type $E_7$, let $A = \{0, 2, 4, 6\}$ and 
$B = \{0, 1, 4, 6\}$.  In this case, $i(A) = b_0 b_2 b_4 b_6$ and
$i(B) = b_0 b_1 b_4 b_6$, $A$ and $B$ are neighbours, and
the equivalence class of $A$ is precisely $\{A, B\}$.
\endexample

\definition{Definition \seca.5 \cite{{\bf 3}, \S6.3}}
Let $n \geq 6$.  

If $n$ is odd, we define 
$P' = P'(n)$ to be the subset of $P(n)$ consisting of the sets $$
\left\{ 
\left\{(n - 1) - 2j : 0 \leq j \leq N \right\}
: 0 \leq N \leq {{n - 1} \over 2} \right\}
,$$ together with the set $$
\{n - 1, n - 3, n - 5, \ldots, 4\} \cup \{0\}
$$ and the empty set.

If $n$ is even, we define 
$P' = P'(n)$ be the subset of $P(n)$ consisting of the sets $$
\left\{
\left\{(n - 1) - 2j : 0 \leq j \leq N \right\}
: 0 \leq N \leq {{n - 2} \over 2} \right\}
,$$ together with the empty set.
\enddefinition

\example{Example \seca.6}
In type $E_6$, we have $$
P' = \left\{
\{5\},
\{5, 3\},
\{5, 3, 1\},
\emptyset
\right\}
.$$
In type $E_7$, we have $$
P' = \left\{
\{6\},
\{6, 4\},
\{6, 4, 2\},
\{6, 4, 2, 0\},
\{6, 4, 0\},
\emptyset
\right\}
.$$ 
\endexample

The importance of the set $P'$ comes from the following

\proclaim{Proposition \seca.7 (Fan, \cite{{\bf 3}, Lemma 8.1.2})}
The set $P'$ constitutes a complete set of equivalence class representatives 
for $P$ with respect to $\sim$.
\qed\endproclaim

\head \secb. Cells and the $\afn$-function \endhead

In \S\secb, we recall the definitions of the $\afn$-function and cells
arising from the monomial basis.  Most of this material comes from the
papers \cite{{\bf 3}} and \cite{{\bf 10}}, or is implicit in them.

\definition{Definition \secb.1 \cite{{\bf 3}, Definition 2.3.1}}
The {\it $\afn$-function} $\afn : W_c \ra \zed^{\geq 0}$
is defined by $$
\afn(w) := \max_{A \in P} \{\#A : w = x i(A) y \text{\ is reduced} \}
$$ for $w \in W_c$.
\enddefinition

\proclaim{Proposition \secb.2}
Let $w \in W_c$ and let $f \in \A$.  Define the 
degree, $\deg f$, of $f$ to be the largest integer $n$ such that $v^n$ occurs
with nonzero coefficient in $f$, with the convention that $\deg 0 = -\infty$.
Denote the structure constants with respect to the monomial basis by 
$g_{x, y, z} \in \A$, namely $$
b_x b_y = \sum_{z \in W_c} g_{x, y, z} b_z
.$$
\item{\rm (i)}
{The structure constant $g_{x, y, z}$ is either zero or a nonnegative power
of $\d$, and, given $x$ and $y$, we have $g_{x, y, z} \ne 0$ for a unique $z$.}
\item{\rm (ii)}
{If $s \in S$ and $g_{s, y, z} \not\in \zed$, then $g_{s, y, z} = \d$,
$\ell(sy) < \ell(y)$ and $y = z$.  Similarly, if  $g_{x, s, z} \not\in 
\zed$, then $g_{x, s, z} = \d$, $\ell(xs) < \ell(x)$ and $x = z$}.
\item{\rm (iii)}
{We have $
\afn(w) = \max_{x, y \in W_c} \deg g_{x, y, w}
.$}
\item{\rm (iv)}
{We have $
\afn(w) = \max_{x, y \in W_c} \deg g_{w, x, y}
.$}
\endproclaim

\demo{Proof}
Parts (i) and (ii) are well known and follow easily from 
\cite{{\bf 3}, Proposition 5.4.1}.

Part (iii) is proved in \cite{{\bf 10}, Proposition 4.2.3} using 
the results of \cite{{\bf 3}}.

The proof of \cite{{\bf 3}, Theorem 5.5.1} shows that $$
\deg g_{w, x, y} \leq \min (\afn(w), \afn(x))
,$$ which means that $$
\max_{x, y \in W_c} \deg g_{w, x, y} \leq \afn(w)
.$$  Conversely, \cite{{\bf 3}, Lemma 5.2.6} shows that $$
b_w b_{w^{-1}} = (v + v^{-1})^{\afn(w)} b_d
$$ for some $d \in W_c$, so taking $x = w^{-1}$ and $y = d$, we find that $$
\max_{x, y \in W_c} \deg g_{w, x, y} \geq \afn(w)
,$$ which completes the proof of (iv).
\qed\enddemo

\definition{Definition \secb.3 \cite{{\bf 3}, Definition 4.1}}

For any $w, w' \in W_c$, we say that $w' \leq_L w$ if there exists $b_x$ 
such that $g_{x, w, w'} \ne 0$, where $g$ is as in Proposition \secb.2.

For any $w, w' \in W_c$, we say that $w' \leq_R w$ if there exists $b_x$ 
such that $g_{w, x, w'} \ne 0$.

For any $w, w' \in W_c$, we say that $w' \leq_{LR} w$ if there exist $b_x$
and $b_y$ such that $b_x b_w b_y = c b_{w'}$ for some $c \ne 0$.

We write $w \sim_L w'$ to mean that both $w' \leq_L w$ and $w \leq_L w'$.
Similarly, we define $w \sim_R w'$ and $w \sim_{LR} w'$.

The relation $\sim_L$ (respectively, $\sim_R$, $\sim_{LR}$) is an equivalence
relation, and the corresponding equivalence classes of $W_c$ are called
the {\it left} (respectively, {\it right}, {\it two-sided}) cells.  

It is clear from the definitions and the fact that the identity element 
is a monomial basis element that two-sided cells are unions of left 
cells, and also unions of right cells.
\enddefinition

\proclaim{Proposition \secb.4}
\item{\rm (i)}
{Let $w \in W_c$.  If we have $w = x i(A) y$ reduced for some $A$ such that
$\#A = \afn(w)$, then $i(A) \sim_{LR} w$ and $w \sim_R x i(A)$.}
\item{\rm (ii)}
{The $\afn$-function is constant on left, right, and two-sided cells.}
\item{\rm (iii)}
{If $w, w' \in W_c$ are such that $w' \leq_R w$ and $w' \not\sim_R w$, then
$\afn(w') > \afn(w)$.  An analogous statement holds for left cells and
two-sided cells.}
\item{\rm (iv)}
{The right cell containing $i(A)$ is precisely the set $$
\{ w \in W_c : w = i(A) x \text{ reduced}, \ \afn(w) = \#A \}
.$$}\item{\rm (v)}
{A left cell and a right cell contained in the same two-sided cell intersect
in a unique element.}
\endproclaim

\demo{Proof}
Statement (i) is proved during the argument establishing \cite{{\bf 3}, Theorem
4.5.1.}.

The fact that the $\afn$-function is constant on two-sided cells is 
implicit in the proof of \cite{{\bf 3}, Theorem 4.5.1}.  Since two-sided cells
are unions of left (or right) cells, part (ii) follows.

Suppose now that $w, w' \in W_c$ are such that $w' \leq_R w$ and 
$w' \not\sim_R w$.  An inductive argument using the definition of $\leq_R$
reduces the problem to the case where there
is some $s \in S$ such that $b_w b_s$ is a multiple of $b_{w'}$, so let
us assume that this is the situation.  By \cite{{\bf 3}, Corollary 4.2.2},
the assumption that $w' \leq_R w$ implies that $\afn(w') \geq \afn(w)$.
The statement follows unless $\afn(w') = \afn(w)$, so suppose we are in this
case.  

Let us write $w = x i(A) y$ as in statement (i).
Now \cite{{\bf 3}, Lemma 4.2.5}, applied to the element $x i(A)$ and the
sequence of generators corresponding to $ys$, shows that we have $w' = 
x i(A) y'$ reduced.  By part (i),
we find that $w' \sim_R x i(A)$, and thus that $w' \sim_R w$,
a contradiction.

The statement for left cells follows by a symmetrical argument, and the
statement for two-sided cells follows from the previous claims and the 
fact that if $w' \leq_{LR} w$, then there is a chain $$
w' = w_1, w_2, w_3, \ldots, w_k = w
$$ where, for each $1 \leq i < k$, we have either $w_i \leq_L w_{i+1}$ or
$w_i \leq_R w_{i+1}$.  This completes the proof of (iii).

Part (iv) is \cite{{\bf 3}, Proposition 4.4.3}.

Part (v) is well known and follows from the proof of \cite{{\bf 3}, Theorem 6.1.2}.
\qed\enddemo

\remark{Remark \secb.5}
For finite and affine Weyl groups, the $\afn$-function defined
above is known by \cite{{\bf 23}, Theorem 3.1} to be
the restriction of Lusztig's more general $\afn$-function 
\cite{{\bf 21}} restricted to the subset $W_c$.

Although it is not true that each of the monomial cells studied above is a cell
in the sense of Kazhdan--Lusztig \cite{{\bf 19}}, it can be shown fairly easily
that each left (respectively, right, two-sided) monomial cell is a subset of
some left (respectively, right, two-sided) Kazhdan--Lusztig cell.
\endremark

\head \secc. Traces on the algebras $TL(E_n)$ \endhead

In \S\secc, we will extend scalars and deal with a $\qvk$-form of
$TL(E_n)$, where $\qvk$ is a field containing $\A$ and a square root of $\d$.
(The existence of $\sqrt{\d}$ is needed for compatibility with \cite{{\bf 3}},
but can ultimately be removed; see Remark \secd.4.)
We write $\tlqvk := \qvk \otimes_\A TL(E_n)$.  We aim to classify the
traces, $\t : \tlqvk \ra \qvk$, that is, linear functions $\t$ with
the property that $\t(ab) = \t(ba)$ for all $a, b \in \tlqvk$.  It is
clear that the set of all traces on $\tlqvk$ is a $\qvk$-vector space
(dependent in principle on $K$ and $\delta$).  The
main result of \S\secc\  is that there is a basis for this vector space
in natural bijection with the set $P'$ of \S\seca.

The next result shows how $\t$ naturally induces a function 
$P/\sim \ra \qvk$.

\proclaim{Lemma \secc.1}
Maintain the notation of Definition \seca.3.  Suppose $A, B \in P$ are such
that $A \sim B$, and let $\t : \tlqvk \ra \qvk$ be a trace.  Then 
$\t(i(A)) = \t(i(B))$.
\endproclaim

\demo{Proof}
The proof immediately reduces to the case where $A$ and $B$ are neighbours.
Let $s$ (respectively, $t$) be the element of $S$ corresponding to the 
unique element of $A \backslash B$ (respectively, $B \backslash A$).  
It is immediate from the 
definitions that $i(A) = b_s i(A \cap B) = i(A \cap B) b_s$ and 
$i(B) = b_t i(A \cap B) = i(A \cap B) b_t$.
We then have $$\eqalign{
\t(i(A)) &= \t(b_s i(A \cap B))
= \t(b_s b_t b_s i(A \cap B)) \cr
&= \t(b_t b_s i(A \cap B) b_s) 
= \t(b_t b_s b_s i(A \cap B)) \cr
&= \d \t(b_t b_s i(A \cap B)) \cr
&= \t(b_t b_t b_s i(A \cap B)) 
= \t(b_t b_s i(A \cap B) b_t) \cr
&= \t(b_t b_s b_t i(A \cap B)) 
= \t(b_t i(A \cap B)) \cr
&= \t(i(B)), \cr
}$$ as required.
\qed\enddemo

\proclaim{Lemma \secc.2}
Any trace $\t : \tlqvk \ra \qvk$ is determined by its values on the set $$
\{ i(A) : A \in P \}
.$$
\endproclaim

\demo{Proof}
Suppose the values of $\t(i(A))$ are known for each $A \in P$.  We will show
how to compute the value of $\t(b_w)$, where $w \in W_c$ is arbitrary.

Let us write $w = x i(A) y$ reduced as in Proposition \secb.4 (i).  
Using a reverse induction, we will assume that the values of $\t(b_{w'})$ for
$\afn(w') > \afn(w) = \#A$, if such $w'$ exist, have been determined.
By the defining relations of $TL(E_n)$, we have $b_{i(A)} b_{i(A)} = 
\d^{\#A} b_{i(A)}$, and so we have $$\eqalign{
\t(b_w) &= \t(b_x b_{i(A)} b_y) \cr
&= \d^{-\#A} \t(b_x b_{i(A)} b_{i(A)} b_y) \cr
&= \d^{-\#A} \t(b_{i(A)} b_y b_x b_{i(A)}). \cr
}$$  Now $i(A) y$ and $x i(A)$ lie in $W_c$ because $w$ does, and Proposition
\secb.4 (i) and (ii) shows that $\afn(w) = \afn(x i(A))$.  By Proposition
\secb.2 (i), we have $$
b_{i(A)} b_y b_x b_{i(A)} = b_{i(A) y} b_{x i(A)} = \d^c b_z
$$ for some $z \in W_c$, and it is clear from the definitions that $z \leq_L 
x i(A)$.  By Proposition \secb.4 (ii) and (iii), we see that $$
\afn(z) \geq \afn(x i(A)) = \afn(w) = \#A
.$$  

If $\afn(z) > \#A$ then our inductive hypothesis determines the value
of $\t(\d^c b_z)$, which in turn determines the value of $\t(b_w)$.  We
may therefore assume that $\afn(z) = \#A$.  To complete the proof, it is
enough to show that $z = i(A)$, because the value of $\t(b_z)$ will then have
been determined by our assumptions.

Let $s \in A$.  Since $b_s b_{i(A)} = \d b_{i(A)}$ by the defining relations,
the definition of $b_z$ shows that $b_s b_z = \d b_z$.  By Proposition 
\secb.2 (ii), this means that $\ell(sz) < \ell(z)$, and it follows that 
$A \subseteq
\L(z)$.  Because $A$ is a set of commuting generators, standard properties
of Coxeter groups show that we can write $z = i(A) z'$ reduced. 
Applying Proposition \secb.4 (iv) to the fact that $\afn(z) = \#A$
shows that $z \sim_R i(A)$.  A symmetrical argument then shows that we have 
$z \sim_L i(A)$.  By Proposition \secb.4 (v), this can only happen if
$z = i(A)$.
\qed\enddemo

\proclaim{Theorem \secc.3}
For each $\bar{A} \in \bar{P}$ (as in Definition \seca.3), there is a unique
trace $\t_{\bar{A}} : \tlqvk \ra \qvk$ such that for each $B \in P$ we have $$
\t_{\bar{A}}(i(B)) = \cases
1 & \text{ if } B \in \bar{A},\cr
0 & \text{ otherwise.}\cr
\endcases
.$$  The set $$\{\t_{\bar{A}} : \bar{A} \in \bar{P}\}$$ is a $\qvk$-basis for 
the set of all traces $\t : \tlqvk \ra \qvk$.
\endproclaim

\demo{Proof}
It is clear from the definition of trace that the traces from $\tlqvk$ to
$\qvk$ form a $\qvk$-vector space.  Lemmas \secc.1 and \secc.2 show that this
space has dimension at most the size of $\bar{P}$.  

Fan \cite{{\bf 3}, Theorem 5.6.1} shows that $\tlqvk$ is semisimple 
and that is then a direct sum of $|\bar{P}|$ matrix rings.  This proves
that the dimension of the space of traces is at least the size of 
$\bar{P}$, and thus that the space has the claimed dimension.

A dimension count, together with another application of lemmas \secc.1 and
\secc.2, then shows that there are unique traces $\t_{\bar{A}}$ with the
properties claimed, and that they form a basis.
\qed\enddemo

We now come to the central definition of the paper.

\definition{Definition \secc.4}
The trace $\tr : \tlqvk \ra \qvk$ is defined by $$
\tr = \sum_{{\bar A} \in {\bar P}} \d^{-\#A} \t_{\bar{A}}
,$$ where $\t_{\bar{A}}$ is as in Theorem \secc.3.
\enddefinition

\proclaim{Corollary \secc.5}
Any trace $\t : \tlqvk \ra \qvk$ satisfies $\t(b_w) = \t(b_{w^{-1}})$
for all $w \in W_c$.
\endproclaim

\demo{Proof}
It follows from Proposition \seca.1 that there is a unique $\A$-linear 
antiautomorphism $* : TL(E_n) \ra TL(E_n)$ fixing the generators $b_s$.
We may extend this to a $\qvk$-linear antiautomorphism $* : \tlqvk 
\ra \tlqvk.$  If $a \in \tlqvk$, let us write $a^*$ for $*(a)$.  Note that
if $A \in P$, then $i(A)$ is invariant under $*$, because $i(A)$ is a product
of commuting generators $b_s$.

Given a trace $\t : \tlqvk \ra \qvk$, the $\qvk$-linear map 
$\t' : \tlqvk \ra \qvk$ defined by $\t'(a) = \t(a^*)$ is also a trace.
Since $\t$ and $\t'$ agree on all elements $i(A)$ for $A \in P$, Lemma
\secc.2 shows that $\t = \t'$, and the assertion follows.
\qed\enddemo

\remark{Remark \secc.6}
The trace $\tr$ will turn out to induce the Markov trace of the title.
Note that the definition makes sense because $\bar{A}, \bar{B} \in \bar{P}$
implies $\#A = \#B$.

Traces on Hecke algebras of finite Coxeter groups are known have a property
similar to that given in Corollary \secc.5; see \cite{{\bf 5}, Corollary 8.2.6}
for more details.
\endremark

\head \secd. Cellular structure and the $\afn$-funtion \endhead

In \S\secd, we explain how the trace $\tr$ is particularly compatible with
the structure of $TL(E_n)$ as a cellular algebra, in the sense of \cite{{\bf 8}}.
We will not recall the complete definition of a cellular algebra here,
but we summarize below the properties of the cellular structure that are
important for our purposes.

\definition{Definition \secd.1}
Let $\Lambda$ be the set of two-sided cells for $TL(E_n)$, equipped with
the partial order induced by $\leq_{LR}$.  For each $\l \in \Lambda$, let
$M(\l)$ be an indexing set for the left cells contained in $\l$; note that
the inversion map on the Coxeter group $W$ induces a bijection between
the set of left cells in $\l$ and the set of right cells in $\l$ (see
the remarks at the end of \cite{{\bf 3}, \S4.4}).
\enddefinition

\proclaim{Proposition \secd.2}
Maintain the above notation.
\item{\rm (i)}{Let $T, U \in M(\l)$ for some fixed $\l \in \Lambda$.  Then
$T \cap U$ contains a unique element, $w$, and we define $C_{T, U} = b_w$.}
\item{\rm (ii)}{The $\A$-algebra anti-automorphism $* : TL(E_n) \ra TL(E_n)$
defined by $*(b_w) = b_{w^{-1}}$ satisfies $*(C_{T, U}) = C_{U, T}$.  In 
particular, we have $w^2 = 1$ if and only if $b_w = C_{T, T}$ for some $T$.}
\item{\rm (iii)}{Suppose that $C_{P, Q}$ and $C_{R, S}$ are arbitrary 
monomial basis elements, and define $C_{T, U}$ by the condition $$
C_{P, Q} C_{R, S} = \d^a C_{T, U}
$$ (which makes sense by Proposition \secb.2 (i)).  If $P, Q, R, S, T$ and 
$U$ all belong to the same two-sided cell, then $P = T$ and $S = U$; if,
furthermore, we have $Q = R$, then $a = \afn(C_{T, U})$.  If it is not the
case that $P = T$, $S = U$ and $Q = R$, then we have $a < \afn(C_{T, U})$.}
\endproclaim

\demo{Proof}
Parts (i) and (ii), which are originally due to Graham \cite{{\bf 7}}, are 
proved in \cite{{\bf 10}, Proposition 4.2.1}.  Part (iii) is proved in
\cite{{\bf 10}, propositions 4.2.1 and 4.2.3} using the results of \cite{{\bf 3}}.
\qed\enddemo

\proclaim{Proposition \secd.3}
For all $w \in W_c$, we have $\tr(b_w) = \d^a$, where $a = -\afn(w)$ if
$w^2 = 1$, and $a < -\afn(w)$ otherwise.
\endproclaim

\demo{Proof}
Let $\l$ be the two-sided cell containing $w$.  We will prove the statement
by induction on the partial order on two-sided cells given in Definition
\secd.1.  Writing $w = C_{T, U}$ for $T, U \in M(\l)$, as in Proposition
\secd.2 (i), and applying Proposition \secd.2 (ii),  we see that the 
condition $w^2 = 1$ is equivalent to $T = U$.

By Proposition \secb.4, there exists a product of $\afn(w)$ 
commuting generators,
$i(A)$, in $\l$.  Define $V \in M(\l)$ by the condition $C_{V, V} = b_{i(A)}$.
Since $\tr$ is a trace, Proposition \secd.2 (iii) shows that $$
\tr(C_{T, U}) = \d^{-\afn(w)} \tr(C_{T, V} C_{V, U})
= \d^{-\afn(w)} \tr(C_{V, U} C_{T, V})
.$$  By Proposition \secb.2 (i), we have $$
C_{V, U} C_{T, V} = \d^b C_{X, Y}
$$ for some $b \geq 0$ and some basis element $C_{X, Y}$.  There are now
two cases to consider.

The first possibility is that $C_{X, Y}$ comes from the two-sided cell $\l$.
(If $T = U$, this case must occur by Proposition \secd.2 (iii).)
In this case, we have $X = Y = V$, and thus $C_{X, Y} = b_{i(A)}$. 
Proposition \secd.2 (iii) then shows that
$b = \afn(w)$ if $T = U$, and $b < \afn(w)$ otherwise.  Since we
have $\tr(b_{i(A)}) = \d^{-\afn(w)}$ by definition of $\tr$, we have 
$\tr(C_{T, U}) = \d^{-\afn(w) + b - \afn(w)}$, and the result follows.

The other possibility is that $C_{X, Y}$ comes from a two-sided cell $\l'$
with $\l' < \l$, and $T \ne U$.  In this case, Proposition \secb.4 (iii) 
shows that $\afn(C_{X, Y}) > \afn(w)$.  By the inductive hypothesis, we 
know that 
$\tr(C_{X, Y}) = \d^{a'}$, where $a' \leq -\afn(C_{X, Y}) < -\afn(w)$.
This means that $\tr(C_{T, U}) = \d^{- \afn(w) + b + a'}$.  By
propositions \secb.2 (iii) and \secb.4 (ii), we have $b \leq \afn(w)$, and 
thus $\tr(C_{T, U}) = \d^a$ for $a < -\afn(w)$, as required.
\qed\enddemo

\remark{Remark \secd.4}
The above proposition shows that we do not actually need $\sqrt{\d} \in k$
to define $\tr$.  From now on, we need only assume that $\qvk$ is a
field containing $\A$.
\endremark

\proclaim{Proposition \secd.5}
If $\qvk$ is the field of fractions of the power series ring $\zed[[v^{-1}]]$,
then $\tr$ is a nondegenerate trace on $\tlqvk$, and $$
\tr(C_{P, Q} C_{R, S}) - \d_{QR} \d_{PS} \in v^{-1} \kyu[[v^{-1}]]
,$$ where $\d_{QR}$ and $\d_{PS}$ are the Kronecker delta.
\endproclaim

\demo{Proof}
An element $x$ of $\qvk$ is uniquely representable in the form $$
x = \sum_{i = -\infty}^N \l_i v^i
,$$ where $\l_i \in \kyu$ for all $i$.  
If $x \ne 0$, we define $\deg x$ to be the largest integer $j$
such that $\l_j \ne 0$.  If $x, y \ne 0$ then $\deg (xy) = \deg x + \deg y$,
so the facts that $\deg \d = 1$ and $\deg 1 = 0$ imply that $\deg \d^a = -a$.

The second assertion follows from 
the fact that $\deg \d^a = -a$ combined with Proposition \secb.4 (ii), 
Proposition \secd.2 (iii) and Proposition \secd.3.

We will now show that for any nonzero $a \in \tlqvk$, we have 
$\tr(a a^*) \ne 0$, from which the assertion follows.  
We have $$
a = \sum_{w \in W_c} \l_w b_w
,$$ and by clearing denominators (thus multiplying $a$ by a nonzero scalar), 
we may assume that we have $\l_w \in \A$ for all $w \in W_c$.  Choose
$w'$ with $\l_{w'} \ne 0$ and $N(w') := \deg \l_{w'}$ maximal, and 
let $c_{w'}$ be the (integer) 
coefficient of $v^{N(w')}$ in $\l_{w'}$.  Setting $a_{w'} = v^{-N(w')} \l_{w'}
b_{w'}$, we then have $$
\tr(a_{w'} a_{w'}^*) = c^2 \mod v^{-1} \kyu[[v^{-1}]]
.$$  
If $\l_{w''} \ne 0$ but $\deg \l_{w''}$ is not maximal, we may again define
$a_{w''} = v^{-N(w'')} \l_{w''} b_{w''}$, but then $$
\tr(a_{w''} a_{w''}^*) \in v^{-1} \kyu[[v^{-1}]]
.$$  Since the integers $c^2$ are strictly positive, it follows that $$
\tr((v^{-N(w')}a)(v^{-N(w')}a)^*) \not\in v^{-1} \kyu[[v^{-1}]]
,$$ which completes the proof.
\qed\enddemo

\proclaim{Proposition \secd.6}
Let $\qvk$ be the field of fractions of the power series ring $\zed[[v^{-1}]]$,
and let $\qvk'$ be the subfield of $\qvk$ consisting of the field of fractions
of $\zed[[v^{-2}]]$.
\item{\rm (i)}{The field $\tlqvk$ has a unique structure as a 
$\zed_2$-graded algebra over $\qvk'$ in which $v^n$ has degree $n \mod 2$ 
and $\qvk'$ is precisely the set of elements of degree $0 \mod 2$.}
\item{\rm (ii)}{The algebra $\tlqvk$ has a unique structure as a 
$\zed_2$-graded
algebra over $\qvk'$ in which $v^n$ has degree $n \mod 2$ and the generators
$b_s$ have degree $1 \mod 2$.  We denote the even subalgebra 
consisting of elements of degree $0 \mod 2$ by $\tlqvkd$.}
\item{\rm (iii)}{Let $\t : \tlqvk \ra \qvk$ be any trace.  Then there are
unique $\qvk'$-linear maps $\t_{(0)}, \t_{(1)} : \tlqvkd \ra \qvk'$ such
that $\t_{(0)} + v \t_{(1)}$ is the restriction of $\t$ to $\tlqvkd$, and
furthermore, $\t_{(0)}$ and $\t_{(1)}$ are themselves traces.}
\endproclaim

\demo{Proof}
Recall from the proof of Proposition \secd.5 that
$\qvk = \kyu((v^{-1})) = \kyu[v][[v^{-1}]]$, so that each element $x \in \qvk$
has a unique expression of the form $$
\sum_{i = N}^\infty q_i v^i
,$$ where $q_i \in \kyu$ and $N \in \zed$ depends on $x$.  Similar reasoning
shows that the subfield $\qvk'$ of $\qvk$ then consists precisely of those 
elements for which $q_i = 0$ whenever $i$ is odd.  Part (i) is a consequence
of this construction.

The assertion of (ii) is immediate from the observation that the defining
relations of Proposition \seca.1 respect the given grading.

Let $\pi : \qvk \ra \qvk'$ be the map $$
\pi\left( \sum_{i = N}^\infty q_i v^i\right) =  \sum_{i = N}^\infty 
q'_i v^i
,$$ where $$q'_i = \cases
q_i & \text{ if } i \text{ is even,}\cr
0 & \text{ otherwise.}\cr
\endcases$$  Our description of $\qvk'$ shows that $\pi$ is a $\qvk'$-linear
map.  Denoting the restriction of $\t$ to $\tlqvkd$ by $\t'$, it follows
that $\pi \circ \t'$ is a trace on $\tlqvkd$.  Since $\t_{(0)} = 
\pi \circ \t'$, the maps $\t_{(0)}$, $v\t_{(1)} = \t' - \t_{(0)}$ and 
$\t_{(1)}$ are also traces, completing the proof of (iii).
\qed\enddemo

Note that any trace from $\tlqvkd$ to $\qvk'$ extends uniquely to a trace
from $\tlqvk$ to $\qvk$ by tensoring by $\qvk \otimes_{\qvk'} -$.

\proclaim{Lemma \secd.7}
The trace $\tr : \tlqvk \ra \qvk$ arises from a trace $$\tr' : \tlqvkd \ra
\qvk'$$ by extension of scalars.
\endproclaim

\demo{Proof}
We use the notation of \S\secc.
Note that if $A \in P$, then $i(A)$ is an element of $\tlqvk$ of degree
$\#A \mod 2$.  We also have $\tr(i(A)) = \d^{\#A}$, which is an element
of $\qvk$ of degree $\#A \mod 2$.

Recall that $\tlqvkd$ is a $\qvk'$-subalgebra of $\tlqvk$ and note
that if $y, z$ are homogeneous elements of $\tlqvk$, then $yz$ and $zy$ have
the same degree.
The argument of Lemma \secc.2 now shows that if $x$ is an element of
$\tlqvkd$, we have a relation $$
\tr(x) = \tr\left( \sum_{\bar{A} \in \bar{P}} \l_{\bar{A}} i(A) \right)
,$$ where for each $\bar{A} \in \bar{P}$, we have $\l_{\bar{A}} (i(A)) 
\in \tlqvkd$.  By the first paragraph of the proof, $\l_{\bar{A}}$ must be
homogeneous of degree $\#A \mod 2$, and $\tr(\l_{\bar{A}} i(A)) \in \qvk'$.
The proof is completed by the observation that any $x \in \tlqvk$ is uniquely
expressible as $x_{(0)} + v x_{(1)}$ for $x_{(0)}, x_{(1)} \in \tlqvkd$
(compare with Proposition \secd.6 (iii)).
\qed\enddemo

\proclaim{Corollary \secd.8}
If $w \in W_c$ and $\tr(b_w) = \d^a$ as in Proposition \secd.3, then 
$a \equiv \l(w) \mod 2$.
\endproclaim

\demo{Proof}
By Lemma \secd.7, we have $\deg\tr(b_w) = \ell(w) \mod 2$, so the assertion
follows from the fact that $\deg \d = 1$.
\qed\enddemo

\head \S\sece.  tom Dieck's diagram calculus \endhead

In \cite{{\bf 2}}, tom Dieck introduced a diagram calculus for the algebras
$TL(E_n)$.  To give a rigorous definition of tom Dieck's diagram calculus,
as we do here, we first need to recall the graphical definition of the
Temperley--Lieb algebra. We start by recalling Jones' formalism of 
$k$-boxes \cite{{\bf 18}}, following the approach of Martin and the author in
\cite{{\bf 15}}.  For further details 
and references, the reader is referred to \cite{{\bf 11}, \S2}.

\definition{Definition \sece.1}
Let $k$ be a nonnegative integer.  The {\it standard $k$-box}, ${\Cal
B}_k$, is the set $\{(x, y) \in \real^2 : 0 \leq x \leq k + 1, \ 0
\leq y \leq 1\}$, together with the $2k$ marked points $$\eqalign{
&1 = (1, 1), \ 
2 = (2, 1), \ 
3 = (3, 1), \ 
\ldots, \ k = (k, 1), \cr
&k + 1 = (k, 0), \
k + 2 = (k-1, 0), \
\ldots, \ 
2k = (1, 0).\cr
}$$
\enddefinition

\definition{Definition \sece.2}
Let $X$ and $Y$ be embeddings of some topological spaces (such as lines) 
into the standard $k$-box. 
Multiplication of such embeddings to obtain a new embedding in the
standard $k$-box shall, where appropriate, be defined via the following
procedure on $k$-boxes.  
The product $XY$ is the embedding
obtained by placing $X$ on top of $Y$ 
(that is, $X$ is first shifted in the plane by $(0,1)$ relative to $Y$, 
so that marked point $(i,0)$ in $X$ coincides with $(i,1)$ in $Y$),
rescaling vertically by a scalar factor of $1/2$ and applying the
appropriate translation to recover a standard $k$-box.
\enddefinition

\definition{Definition \sece.3}
Let $k$ be a nonnegative integer. Consider the set of 
smooth embeddings of a single curve 
(which we usually call an ``edge'') in the standard $k$-box,
such that the curve is either closed (isotopic to a circle) 
or its endpoints coincide with two marked points of the box,
with the curve meeting the boundary of the box only at such
points, and there transversely.

By a smooth diffeomorphism of this curve we mean a smooth
diffeomorphism of the copy of $\real^2$ in which it is embedded,
that fixes the boundary, and in particular the marked points, of the
$k$-box, and takes the curve to another such smooth embedding. 
(Thus, the orbit of smooth diffeomorphisms
of one embedding contains all embeddings with the same endpoints.)

A concrete Brauer diagram is a set of such embedded curves with the property
that every marked point coincides with an endpoint of precisely one
curve. 
(In examples we can represent this set by drawing all the curves on
one copy of the $k$-box. Examples can always be chosen in which no 
ambiguity arises thereby.)

Two such concrete diagrams are said to be equivalent if one may be
taken into the other 
by applying smooth diffeomorphisms to the individual curve
embeddings within it.

There is an obvious map from the set of concrete diagrams to the set
of pair partitions of the $2k$ marked points. 
It will be evident that the image under this map is an invariant of
concrete diagram equivalence.  

The set $B_k(\emptyset)$ is the set of equivalence classes of
concrete diagrams. Such a class (or any representative) is called
a Brauer diagram.  

Let $D_1, D_2$ be concrete diagrams. 
Since the $k$-box multiplication defined above
internalises marked points in coincident pairs, 
corresponding curve endpoints in $D_1D_2$ may also be internalised seamlessly. 
Each chain of curves concatenated in this way may thus be put in
natural correspondence with a single curve. 
Thus the multiplication 
gives rise to a closed associative binary operation on the
set of concrete diagrams.
It will be evident that this passes to a well defined
multiplication on $B_k(\emptyset)$. 
Let $R$ be a commutative ring with $1$.  
The elements of $B_n(\emptyset)$ form the basis elements of an $R$-algebra
$\Pl_n^B(\emptyset)$ with this multiplication. 
\enddefinition

A curve in a diagram that is not a closed loop
is called {\it propagating} if its endpoints have different
$y$-values, and {\it non-propagating} otherwise.  (Some authors use the
terms ``through strings'' and ``arcs'' respectively for curves of these
types.)

Note that in a Brauer diagram drawn on a single copy of the $k$-box it
is not generally possible to keep the embedded curves disjoint.
Let  $T_k(\emptyset) \subset B_k(\emptyset)$
denote the subset of diagrams having representative elements in which
the curves are disjoint. 
Representatives of this kind are called Temperley--Lieb diagrams.

It will be evident that $\Pl_n^B(\emptyset)$ has a subalgebra with
basis the subset $T_k(\emptyset) $. 
(That is to say, the disjointness property is preserved under multiplication.)
We denote this subalgebra 
$\Pl_n(\emptyset)$ 

Because of the disjointness property there is, for each element of 
$T_k(\emptyset) $, a unique assignment of orientation to its 
curves that satisfies the following two conditions.

\item{(i)}{A curve meeting the $r$-th marked point of the standard
$k$-box, where $r$ is odd, must exit the box at that point.}

\item{(ii)}{Each connected component of the complement of the union of
the curves in the standard $k$-box may be oriented in such a way that
the orientation of a curve coincides with the orientation induced as
part of the boundary of the connected component.}

Note that the orientations match up automatically in composition.  If
$D_1$ and $D_2$ are equivalent concrete Temperley--Lieb diagrams, the
diffeomorphisms that give rise to the equivalence set up a bijection
between the connected components of $D_1$ and those of $D_2$.

\topcaption{Figure 2} A pillar diagram corresponding to an element of
$T_8(\emptyset)$ \endcaption
\centerline{
\hbox to 4.6in{
\vbox to 2.2in{\vfill
        \includegraphics{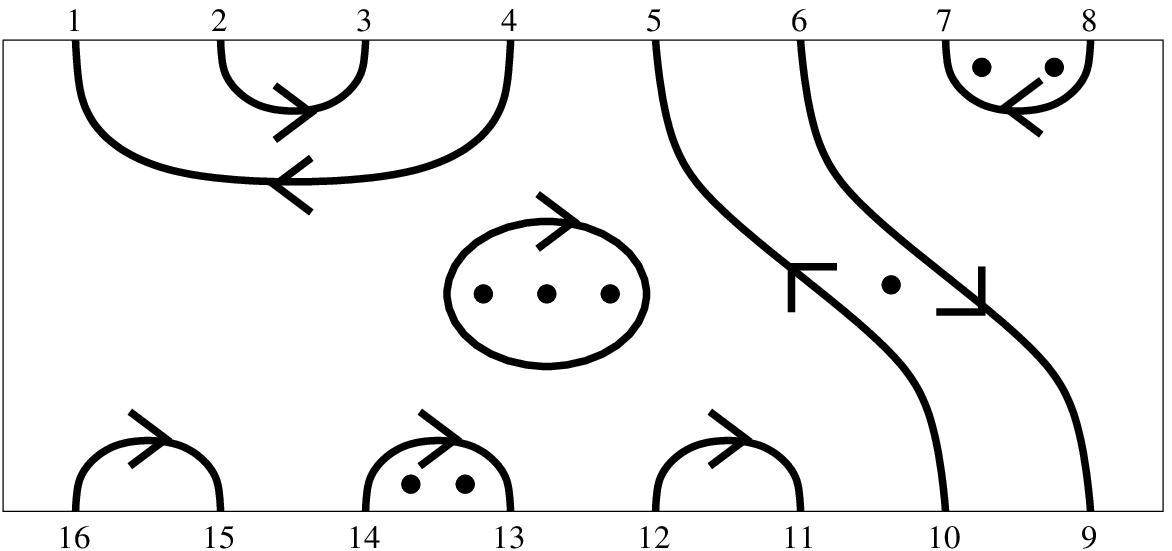}
}
\hfill}
}

\definition{Definition \sece.4}
A {\it pillar diagram} consists of a pair $(D, f)$, where 
$D \in T_k(\emptyset)$ is a Temperley--Lieb diagram and $f$ is 
a function from the connected components of $D$ to $\zed^{\geq 0}$, such
that any component with anticlockwise orientation is mapped to zero.

On the diagram $D$, we indicate the values of $f$ on the clockwise
connected components either by writing in the appropriate integer, or by
inserting $k$ disjoint discs (the ``pillars'' of \cite{{\bf 2}}).

The set of pillar diagrams arising from the set $T_k(\emptyset)$ will be
denoted $T_k(\bullet)$.
\enddefinition

\example{Example \sece.5}
Let $k = 8$.  A pillar diagram corresponding to an element of
$T_k(\bullet)$ is shown in Figure 2.
Note that there are 10 connected components,
precisely 7 of which inherit a clockwise orientation.  The values of $f$
on these 7 components are $3, 2, 2, 1, 0, 0, 0$.
\endexample

We define an algebra $\Pl_n(\bullet)$, analogous to $\Pl_n(\emptyset)$,
with the set $T_k(\bullet)$ as a basis.  The multiplication is $k$-box
multiplication with the added convention that function values on the
connected components are additive.  (This is natural if one represents
the function values with pillars as in Figure 2.)

For our purposes, we need to apply an equivalence relation on the concrete
diagrams of $T_k(\bullet)$.  Locally, this is given by the relation shown
in Figure 3.

\topcaption{Figure 3} A topological reduction rule \endcaption
\centerline{
\hbox to 2.263in{
\vbox to 1.000in{\vfill
        \includegraphics{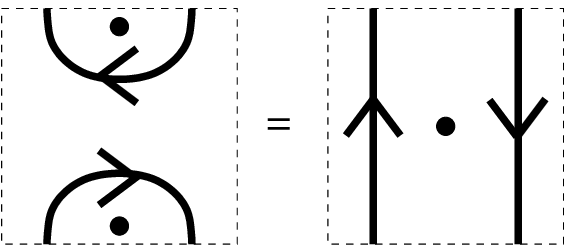}
}
\hfill}
}

In the notation where clockwise regions are labelled by nonnegative integers,
the relation of Figure 3 is that shown in Figure 4.

\topcaption{Figure 4} Alternative notation for the topological reduction 
\endcaption
\centerline{
\hbox to 2.263in{
\vbox to 1.000in{\vfill
        \includegraphics{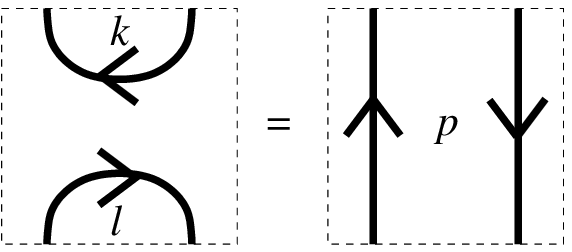}
}
\hfill}
}

If the regions labelled $k$ and $l$ are connected to each other, Figure 3
shows that we have $k = l > 1$ and $p = k-1$.  On the other hand, if the 
regions labelled $k$ and $l$ are genuinely distinct, that is, the arcs shown
on the left hand side of figure 3 are not sections of some longer arc,
then we have $p = k + l - 1 \geq 1$.  In the latter case, it is not possible 
for any regions labelled by the integer zero to be created or destroyed by
the topological reduction.
Note that the other partial regions shown in figures 2 and 3 have 
anticlockwise orientation, and as such they are labelled by the integer $0$.

\definition{Definition \sece.6}
If $L$ is a closed loop in a concrete diagram of $T_k(\bullet)$, we define
$m(L)$ to be the integer label of the region immediately interior to $L$;
in particular, we have $m(L) = 0$ if $L$ has anticlockwise orientation.

Let $R$ be a commutative ring with $1$.  The $R$-algebra $\Pl_n^E(\bullet)$
is the quotient of the $R$-algebra $\Pl_n(\bullet)$ obtained by applying
the following three relations:
\item{(i)}{for each closed loop $L$ whose immediate interior is labelled
$1$ and whose immediate exterior is necessarily labelled $0$, relabel the 
immediate interior of $L$ by $0$ and remove $L$;}
\item{(ii)}{for each closed loop $L$ whose immediate interior is labelled
$0$ and whose immediate exterior is labelled $k$, relabel the immediate
interior of $L$ by $k$, remove $L$ and multiply by $\d$;}
\item{(iii)}{for each region $R$ labelled by $k \geq 2$ (whether or not $R$ 
is a closed loop), decrease the label of $R$ by $1$ and multiply by $\d$.}
\enddefinition

A basis for $\Pl_n^E(\bullet)$ may be obtained by using the notion of
``reduced'' diagrams given in \cite{{\bf 2}, \S2} and Bergman's diamond lemma
\cite{{\bf 1}}.  However, we do not pursue this because we do not need it for
our purposes.

\definition{Definition \sece.7}
Suppose $n > 1$ and $1 \leq k < n$.  

The diagram $E_k^n$ of $\Pl_n^E(\bullet)$ is the one where each
point $i$ is connected by a propagating edge to point $2n + 1 - i$,
unless $i \in \{k, k+1, 2n - k, 2n + 1 - k\}$.  Points $k$ and $k+1$
are connected by an edge, as are points $2n - k$ and $2n + 1 - k$.  All 
regions are labelled by $0$.

The diagram $B_k^n$ of $\Pl_n^E(\bullet)$ is the one where each
point $i$ is connected by a propagating edge to point $2n + 1 - i$, and
all regions are labelled by $0$, except the rectangular region bounded by
$k, k+1, 2n - k$ and $2n + 1 - k$, which is labelled by $1$.
\enddefinition

\proclaim{Proposition \sece.8}
There is a unique homomorphism $\rho : TL(E_n) \ra \Pl_n^E(\bullet)$ of 
unital $\A$-algebras sending $b_0$ to $B_3^n$ and $b_s$ to $E_s^n$ for
$i \in \{1, 2, \ldots, n-1\}$, where the numbering of generators is as
in \S\seca.
\endproclaim

\demo{Proof}
This is a routine (but important) exercise using the presentation of
Proposition \seca.1, and is essentially the same
as the proof of \cite{{\bf 2}, Theorem 2.5}.
\qed\enddemo

We shall see later that $\rho$ is in fact a faithful representation.
We will not determine the image of $\rho$, but this can be done by
an inductive combinatorial argument similar to those in \cite{{\bf 9}, \S5}.

\head \S\secf.  Existence and uniqueness of the Markov trace \endhead

There is a well-known embedding $\iota_n : TL(E_n) \ra
TL(E_{n+1})$ sending $b_s$ to $b_s$ for each generator of $TL(E_n)$ (see
\cite{{\bf 3}, \S6.3}).  This means that the tower of algebras $TL(E_n)$, equipped
with the generators $b_s$, fits into the framework of Markov traces defined
in \S\secy.  We recall the definition in order to fix some notation.

\definition{Definition \secf.1}
Let $\qvk$ be a field containing $\A$.
A {\it Markov trace} on $TL_K(E_\infty)$ with parameter $z \in \qvk$
is a $\qvk$-linear map 
$\t : TL_K(E_\infty) \ra K$ 
satisfying the following conditions:
\item{\rm (i)}{$\t(1) = 1$;}
\item{\rm (ii)}{$\t(h b_n) = z \t(h)$ for $n \geq 6$ and $h \in \tlqvk$;}
\item{\rm (iii)}{$\t(hh') = \t(h'h)$ for all $n \geq 6$ and $h, h' \in
\tlqvk$.}
\enddefinition

\remark{Remark \secf.2}
Note that in condition (ii), $b_n$ is the unique generator in $TL(E_{n+1})$
that does not lie in $TL(E_n)$.  As mentioned in \cite{{\bf 3}, \S2.2}, the 
algebras $TL(E_n)$ are quotients of the Hecke algebras of the Coxeter
groups $W(E_n)$, and $b_s = q^{-1/2}(T_s + 1)$, where the $T_s$ are the
usual generators for the Hecke algebra as given in \cite{{\bf 16}, \S7}.  This
means that the Markov trace can also be regarded as a trace on a tower of
Hecke algebras.
\endremark

\proclaim{Proposition \secf.3}
If $\t$ is a Markov trace on $TL_K(E_\infty)$, then the parameter $z$ must
be equal to $\d^{-1}$, and $\t$ is unique.  Restricted to $TL(E_n)$,
such a Markov trace must agree with the trace $\tr$.
\endproclaim

\demo{Proof}
Let $n \geq 6$.  Part (ii) of Definition \secf.1 shows that $\t(b_{n-1} b_n)
= z \t(b_{n-1})$.  On the other hand, the defining relations and part (iii) 
of the definition show that $$
\t(b_{n-1} b_n) = \d^{-1} \t(b_{n-1} (b_{n-1} b_n))
= \d^{-1} \t(b_{n-1} b_n b_{n-1}) = \d^{-1} \t(b_{n-1})
,$$ proving the assertion about the parameter.

To prove the other assertions, it suffices to show that, regarding $\tlqvk$
as a subalgebra of $TL_K(E_\infty)$, we have $\t(i(A)) = \d^{-\#A}$ for
$A \in P = P(n)$.  Choose such an $A$.  It follows from Definition \seca.3
that for sufficiently large $N \geq n$, and identifying $A$ in the obvious
way with an element of $P(N)$, we can find $B \in P(N)$ with $A \sim B$
and $B \cap \{b_0, b_1, \ldots, b_5\} = \emptyset$.  The first assertion
together with repeated applications of part (ii) of Definition \secf.1 
(and one application of part (i)) now show that $\t(i(B)) = \d^{-\#B}
= \d^{-\#A}$, and Lemma \secc.1 completes the proof.
\qed\enddemo

To prove that the Markov trace on $TL_K(E_\infty)$ exists, we make use of
the diagram calculus, as hinted in \cite{{\bf 2}, \S6}.

\definition{Definition \secf.4}
Let $k$ be a nonnegative integer.  The {\it standard $k$-cone} is obtained
from the standard $k$-box by identifying each pair of points 
$\{(x, 0), (x, 1)\}$ for each $0 \leq x \leq k+1$, and identifying all the 
points in the set $\{(k + 1, y) : 0 \leq y \leq 1\}$.  The standard $k$-cone
is homeomorphic to a closed disc.

Let $D$ be a diagram in $\Pl_k^E(\bullet)$.  The {\it trace diagram}, 
$\overline{D}$, of
$D$ is obtained by identifying the boundary points of the $k$-box bounding
$D$ to form the standard $k$-cone.
\enddefinition

\topcaption{Figure 5} The trace diagram of the pillar diagram in Figure 2
\endcaption
\centerline{
\hbox to 2.527in{
\vbox to 2.083in{\vfill
        \includegraphics{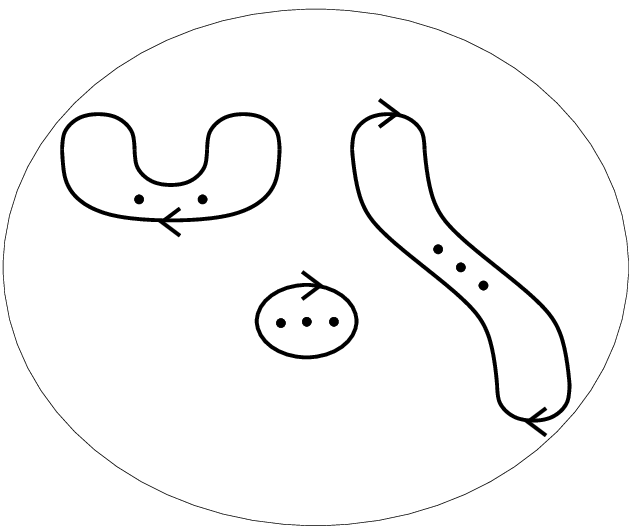}
}
\hfill}
}

\example{Example \secf.5}
The trace diagram $\overline{D}$ 
corresponding to the diagram $D$ of Figure 2 is shown in Figure 5.
\endexample

Notice that the outer part of the trace diagram (regarded as a disc) will
always have an anticlockwise orientation and thus be labelled by $0$.
Consequently, any regions in the trace diagram not labelled by zero 
must be bounded by at least one closed loop.  (It is possible for
the closed loops to be nested.)

\definition{Definition \secf.6}
Let $g : \zed^{\geq 0} \ra \zed^{\geq 0}$ be given by $$
g(c) = \cases
1 & \text{ if } c = 0,\cr
c-1 & \text{ if } c \geq 1.\cr
\endcases
$$  If $\overline{D}$ is a trace diagram for $TL(E_n)$, we define the 
{\it content}, $c(\overline{D})$, of $\overline{D}$ to be the integer $$
\sum_L g(f(L))
,$$ where the sum is over all the connected components $L$ of $\overline{D}$ 
that are interior to at least one closed loop, and where
$f(L)$ is the integer assigned to $L$ as in Definition \sece.4.
\enddefinition

\example{Example \secf.7}
The content of the trace diagram in Figure 5 is $$
g(2) + g(3) + g(3) = 5
.$$
\endexample

\proclaim{Lemma \secf.8}
The content of a trace diagram $\overline{D}$ 
is invariant under the topological reduction rule shown in Figure 3.
\endproclaim

\demo{Proof}
Consider the application of the topological reduction rule to a diagram
that looks locally like the situation in Figure 6.

\topcaption{Figure 6} Labelling of points involved in the topological relation
\endcaption
\centerline{
\hbox to 0.958in{
\vbox to 1.222in{\vfill
        \includegraphics{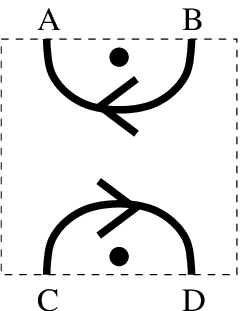}
}
\hfill}
}

As in the discussion following Figure 4, there are two cases to consider, 
according as the two pillar regions are connected or not in $\overline{D}$.

There are four cases to consider, according as there is an oriented curve
in $\overline{D}$ from point A to point C, and (independently) 
according as there is an oriented curve in $\overline{D}$ from point D to
point B.

Suppose first that there is no oriented curve in $\overline{D}$ from point
A to point C, and also that there is no oriented curve in $\overline{D}$
from point D to point B.  In this case, the two pillar regions are
genuinely distinct, and applying the topological relation does not produce
any new closed loops.  We are then in the case $p = k + l - 1 \geq 1$ of
Figure 4, so the summands $(k - 1)$ and $(l - 1)$ appearing in Definition
\secf.6 are replaced by a single $((k + l - 1) - 1)$, leaving the content
unchanged.

We next deal 
with the case where there is an oriented curve from point A to point C,
but no oriented curve from point D to point B.  In this case, the two
pillar regions are connected to each other, and the application of the
topological rule produces a new closed loop (labelled zero) from the curve
originally connecting point A to point C.  We are now in the case
$k = l > 1$ of Figure 4.  This will change one of the summands $(k - 1)$
of Definition \secf.6 to $(k - 2)$, and a new summand of $1$ will be produced,
corresponding to the new closed loop.  The content thus remains unchanged.

Consideration of the case where there is an oriented curve from point D to
point B, but not from point A to point C, proceeds in exactly the same
way.  The last case, in which both oriented curves exist, also works
similarly, except that the oriented curves shown in Figure 6 are already
part of a closed loop.  Application of the topological relation splits this
closed loop into two closed loops, again producing an extra summand of $1$
and changing a summand $(k - 1)$ to $(k - 2)$, leaving the content unchanged.
\qed\enddemo

\proclaim{Lemma \secf.9}
There is a well-defined $K$-linear map $$
\t_n^\bullet : \Pl_n^E(\bullet) \ra K
$$ such that for each pillar diagram $D$, $\t_n^\bullet(D) = 
\d^{c(\overline{D})}$.  If $x, y \in \Pl_n^E(\bullet)$, we have 
$\t_n^\bullet(xy) = \t_n^\bullet(yx)$.
\endproclaim

\demo{Proof}
For the first assertion, we need to check relations (a)--(c) of Definition
\sece.6.  Relation (iii) holds by Lemma \secf.8.

In relation (i), we have $D = D_1$, where $D_1$ is the result of
removing a loop labelled $1$ from $D$.  
Since $c(\overline{D}) = c(\overline{D_1})$, we have 
$\t_n^\bullet(D) = \t_n^\bullet(D_1)$.

In relation (ii), we have $D = \d D_2$, where $D_2$ is the result of
removing a loop labelled $0$ from $D$.  
Since $c(\overline{D}) = c(\overline{D_2}) + 1$, we have 
$\t_n^\bullet(D) = \t_n^\bullet(D_2)$.

By linearity, we only need check the second assertion in the case where $x$
and $y$ are pillar diagrams, and this is immediate from the construction
of trace diagrams from pillar diagrams.
\qed\enddemo

It is not hard to see that there is an algebra embedding $\iota_n^\bullet : 
\Pl_n^E(\bullet) \ra \Pl_{n+1}^E(\bullet)$ analogous to the map $\iota_n$.
Given a pillar diagram $D$ of $\Pl_n^E(\bullet)$, $\iota^\bullet(D)$ is the
diagram obtained by adding a vertical line on the right of the diagram.

\proclaim{Lemma \secf.10}
Let $D$ be a pillar diagram of $\Pl_n^E(\bullet)$.
\item{\rm (i)}{We have $\t_{n+1}^\bullet(\iota_n^\bullet(D)) = 
\d\t_n^\bullet(D).$}
\item{\rm (ii)}{Let $E_n^{n+1}$ be as in Definition \sece.7.  Then we have $
\t_n^\bullet(D) = 
\t_{n+1}^\bullet(\iota_n(D) E_n).
$}
\endproclaim

\demo{Proof}
Part (i) follows from the observation that the trace diagram 
$\overline{\iota_n^\bullet(D)}$ differs from the trace diagram $\overline{D}$
only in having a single extra closed loop, labelled $0$.

A short calculation involving diagrams shows that the trace diagrams
$\overline{D}$ and $\overline{\iota_n(D) E_n}$ are equivalent, from which 
part (ii) follows.
\qed\enddemo

\proclaim{Theorem \secf.11}
Let $\t_n : \tlqvk \ra K$ be the trace defined by $$\t_n(x) = \d^{-n}
\t_n^\bullet(\rho(x)).$$  The family of traces $\{ \t_n : n \geq 6\}$ is
compatible with the direct limit of algebras $\tlqvk$ and gives the unique
Markov trace on $TL_{\qvk}(E_\infty)$.  Furthermore, the Markov trace
agrees with the traces $\tr$ of Definition \secc.4.
\endproclaim

\demo{Proof}
The maps $\t_n$ are traces by Proposition \sece.8 and Lemma \secf.9.  They are
compatible with the direct limit by Lemma \secf.10 (i).  Since $\t_n^\bullet(1)
= \d^n$, we have $\t_n(1) = 1$.  Condition (ii) of Definition \secf.1 follows
from part (ii) of Lemma \secf.10.
Uniqueness of the Markov trace, and agreement with the traces $\tr$,
is given by Proposition \secf.3.
\qed\enddemo

\head \secg. Proofs and applications \endhead

\demo{Proof of Theorem \secz.1}
We need to show that the homomorphism $\rho$ of Proposition \sece.8 is 
injective, and there is no loss in passing to the field of fractions
$\qvk$ of $\zed[[v^{-1}]]$.  In this case, Proposition \secd.5 and Theorem
\secf.11 show that the unique Markov trace on $\tlqvk$, which can be defined
on $\Im(\rho)$, is nondegenerate on $\tlqvk$.  The conclusion follows.
\qed\enddemo

\proclaim{Proposition \secg.1}
The linear map $$
\almd^n \t_n = v^{-n} \t_n^{\bullet} \circ \rho
$$ restricted to $TL(E_n)$ takes values in $\A$.  It is a tabular trace in
the sense of \cite{{\bf 10}}, and a positive generalized Jones trace in the sense of
\cite{{\bf 12}}.
\endproclaim

\demo{Proof}
The first assertion comes from the fact that $\t_n^\bullet$
evaluated on a diagram (such as an element of the form $\rho(b_w)$ for
$w \in W_c$) yields a nonnegative integer power of $\d$.

To check that $\almd^n \t_n$ is a tabular trace, we need to check that 
axiom (A5)
of \cite{{\bf 10}, Definition 1.3.4} is satisfied.  We have just shown that
$\almd^n \t_n$ takes values in $\A$, and it is clear from Theorem \secf.11 that
$\almd^n \t_n$ is a trace.  We have seen in Corollary \secc.5 and 
Proposition \secd.2 (ii) that $\almd^n \t_n(x) = \almd^n \t_n(x^*)$ for all
$x \in TL(E_n)$.  All that remains to check is that $$
\t(v^{\afn(C_{S, T})} C_{S, T}) = \d_{S, T} \mod v^{-1} \A^-
.$$  This follows from propositions \secd.2 (ii) and \secd.3 once we observe
that we have $$
\almd^n = 1 \mod v^{-2} \kyu[[v^{-1}]]
,$$ regarded as power series in $\kyu[v][[v^{-1}]]$.

To show that $\almd^n \t_n$ is a generalized Jones trace (see \cite{{\bf 12}, 
Definition 2.9}), two further conditions must be checked.  One of these
is precisely that established by Lemma \secd.7; the other is that, 
for $x, y \in W_c$, we should have $$
\almd^n \t_n(c_x c_{y^{-1}}) = 
\cases 1 \mod v^{-1} \A^- & \text{ if } x = y,\cr
0 \mod v^{-1} \A^- & \text{ otherwise,} \cr
\endcases$$ where 
$\{c_w : w \in W_c\}$ is the canonical basis of $TL(E_n)$ defined
by J. Losonczy and the author in \cite{{\bf 14}}.  By \cite{{\bf 14}, Theorem 3.6}, 
this is
nothing other than the basis $\{b_w : w \in W_c\}$ in this case.  The
corresponding property for $\tr$ (instead of $\almd^n \t_n$) follows from
Proposition \secd.5, and the assertion for $\almd^n \t_n$ follows from
the fact that $\almd^n = 1 \mod v^{-2} \A^-$.

A generalized Jones trace is positive if it sends canonical basis elements
to elements of $\enn[v, v^{-1}]$.  This holds for $\almd^n \t_n$ by Proposition
\secd.3: in this case, $\almd^n \t_n(b_w) = \d^b$ for some $b \geq 0$, so
that $\almd^n \in \enn[v, v^{-1}]$.
\qed\enddemo

\remark{Remark \secg.2}
Proposition \secg.1 corrects the proof of \cite{{\bf 10}, Theorem 4.3.5}, where 
the proof that the tabular trace takes the same values on $x$ and $x^*$
contains a gap.
\endremark

\demo{Proof of Theorem \secz.2}
By \cite{{\bf 12}, Theorem 7.10}, the conclusion of Theorem \secz.2 holds for
a generalized Jones trace if the underlying Coxeter group has ``Property F''
and a bipartite Coxeter graph.  Clearly the graphs $E_n$ are bipartite,
because they contain no circuits.  Property F holds by \cite{{\bf 12}, Remark 3.5};
see \cite{{\bf 13}, Lemma 5.6} for a fuller explanation.

To complete the proof, we simply have to transfer the result from 
$\almd^n \t_n$ to the Markov trace, which follows from the fact that 
$\almd^n = 1 \mod v^{-2} \A^-$.
\qed\enddemo

The next result is an easier to use version of Theorem \secz.2.

\proclaim{Corollary \secg.3}
Let $x, y \in W_c(E_n)$.  Then we have $$\tmu(x, y) = \cases
1 & \text{ if } \t_n^\bullet \circ \rho(b_x b_{y^{-1}}) = \d^{n-1},\cr
0 & \text{ otherwise.}\cr
\endcases$$
\endproclaim

\demo{Proof}
This follows from Theorem \secz.2 together with 
the observation that $b_x b_{y^{-1}} = \d^b b_w$ for
some $b \geq 0$ and $w \in W_c$, and the fact that $\t_n^\bullet$ sends
diagrams to positive powers of $\d$.
\qed\enddemo

\remark{Remark \secg.4}
It follows from \cite{{\bf 13}, Theorem 4.6 (iv)} and \cite{{\bf 14}, Theorem 3.6} 
that the monomial basis element $b_x$ is the projection of the 
Kazhdan--Lusztig basis element $C'_x \in \H(E_n)$.
Regarding $\tr$ and $\t_n^\bullet \circ \rho$ as traces on the Hecke
algebra, Theorem \secz.2 and Corollary \secg.3 can be used to evaluate the
trace on products of certain Kazhdan--Lusztig basis elements, without 
evaluating the product (which would be difficult).  Another noteworthy
property of these results is that they give non-recursive formulae for
certain of the integers $\mu(x, y)$.
\endremark

\remark{Remark \secg.5}
In \cite{{\bf 7}, \S9}, Graham showed that if $x, w \in W_c$ for $TL(E_n)$ then
$\mu(x, y) \in \{0, 1\}$, and also produced a nonrecursive method of finding
all the $x$ with $\mu(x, y) = 1$ for a fixed $y$.  (In \cite{{\bf 7}}, $x$ and
$y$ are said to be ``close'' if $\tmu(x, y) = 1$.)  However, unlike
the results above, this does not give
an efficient way to compute $\mu(x, y)$ when both of $x$ and $y$ are
specified.  Corollary \secg.3 can therefore be regarded as a quick way to
tell if two elements are close or not.
\endremark

\remark{Remark \secg.6}
It is possible to modify Theorem \secz.2 and Corollary \secg.3 so that they
provide a nonrecursive way to test whether two diagrams represent the same
algebra element.  However, we do not pursue this here for reasons of space.
\endremark

\example{Example \secg.7}
Consider the Coxeter system of type $E_n$ with $n = 6$, and 
generators $s_0, \ldots, 
s_5$ as numbered in Figure 1.  Define $y = s_1 s_2 s_4 s_0 s_5$ and $$
w = s_1 s_2 s_3 s_4 s_0 s_3 s_5 s_2 s_4 s_1 s_3 s_2 s_0 s_3 s_4 s_5
;$$ these are both reduced expressions for fully commutative elements.  
The diagrams $\rho(b_y)$ and 
$\rho(b_w)$ are shown in figures 7 and 8 respectively.  To evaluate
$\t_n^\bullet(b_y b_{w^{-1}})$, we invert the diagram for $b_w$, compose
it with $b_y$ and identify boundary points to produce a trace diagram.
The trace diagram so obtained is shown in Figure 9 (up to equivalence),
and by inspection, it has content $1 + 1 + 1 + (3-1) = 5 = n-1$.  It
follows from Corollary \secg.3 that $\mu(y, w) = 1$.
\endexample

\topcaption{Figure 7} The diagram $\rho(b_y)$ of Example \secg.7 \endcaption
\centerline{
\hbox to 3.541in{
\vbox to 1.958in{\vfill
        \includegraphics{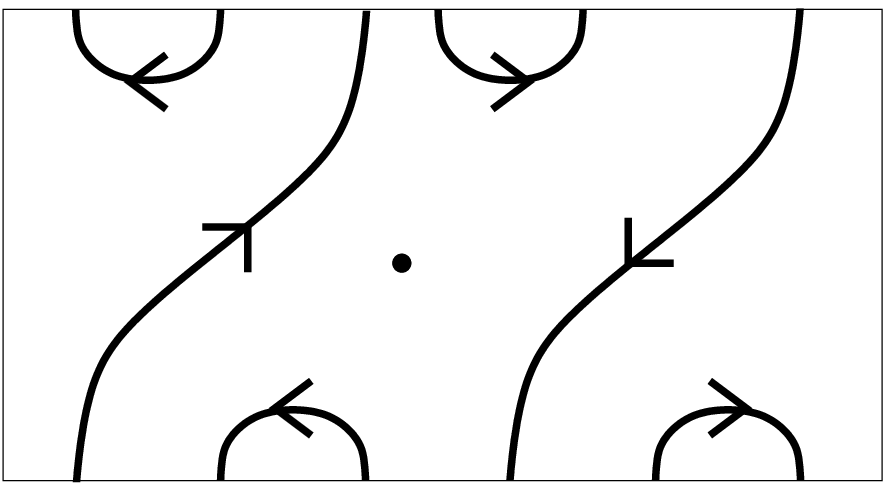}
}
\hfill}
}

\vfill\eject

\topcaption{Figure 8} The diagram $\rho(b_w)$ of Example \secg.7 \endcaption
\centerline{
\hbox to 3.541in{
\vbox to 1.958in{\vfill
        \includegraphics{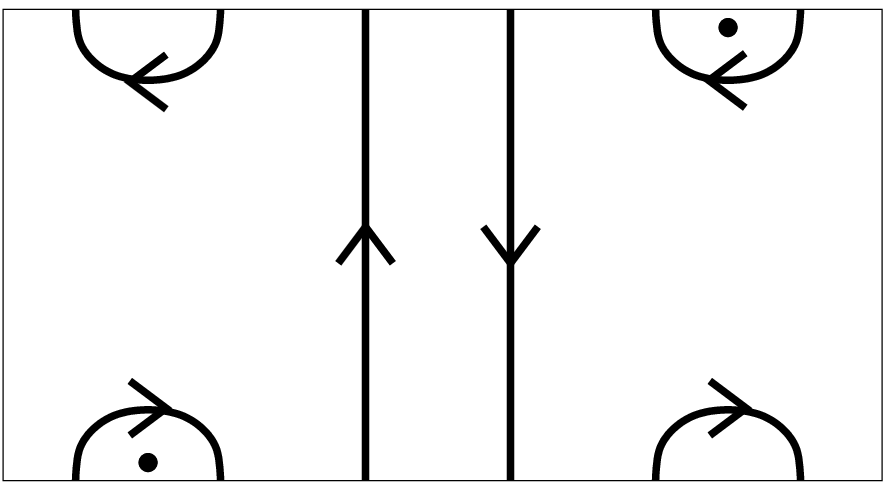}
}
\hfill}
}

\topcaption{Figure 9} The trace diagram corresponding to $\t_6^\bullet \circ
\rho(b_y b_{w^{-1}})$ of Example \secg.7 \endcaption
\centerline{
\hbox to 4.736in{
\vbox to 2.763in{\vfill
        \includegraphics{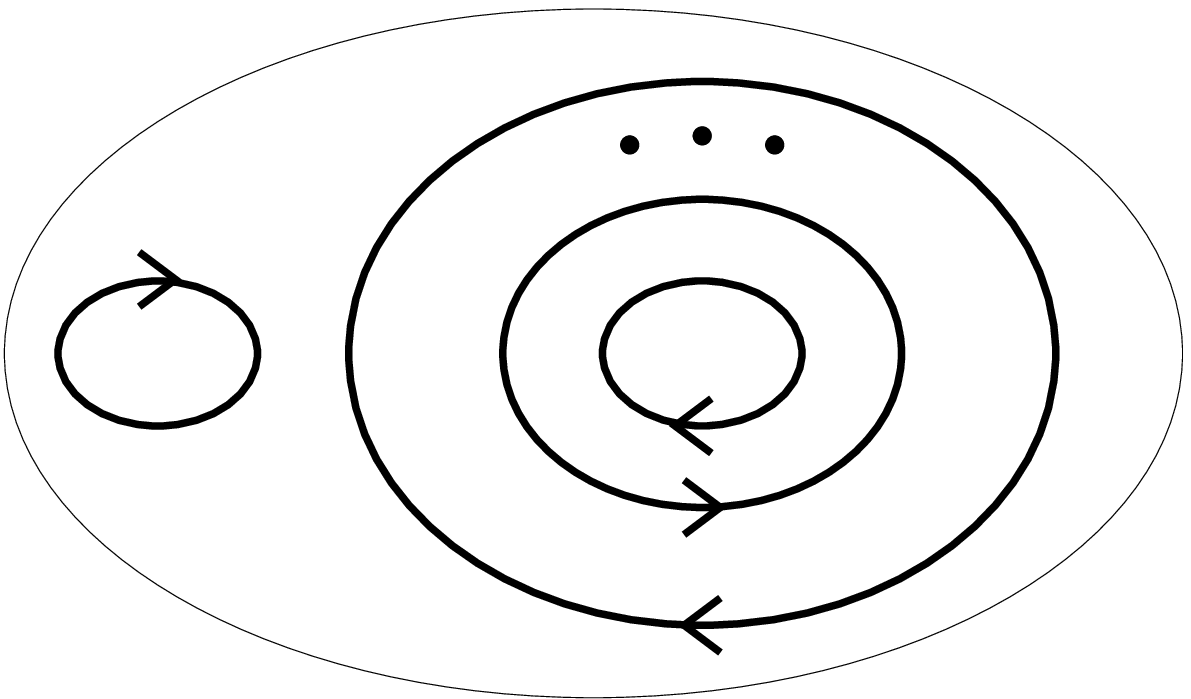}
}
\hfill}
}

\head Acknowledgement \endhead

I am grateful to P.P. Martin for helpful comments on an early version of
this paper.

\leftheadtext{} \rightheadtext{}

\end